\documentclass[12pt]{article}
\usepackage{hyperref}
\usepackage{a4wide}
\usepackage{amsmath}
\usepackage{amssymb}
\usepackage[english,francais]{babel}
\usepackage[fixlanguage]{babelbib}
\usepackage{theorem}
\usepackage{amsfonts}
\usepackage{NN}
\date{}
\title{Sur la d\'efinissabilit\'e existentielle\\ de la non-nullit\'e dans les anneaux\\
{\small \`A para\^{\i}tre dans  \textsl{Algebra and Number Theory} (accept\'e le 19 octobre 2007)}
}
\author{Laurent Moret-Bailly
\thanks{L'auteur est membre du r\'eseau europ\'een `Arithmetic 
Algebraic Geometry' (contrat HPRN-CT-2000-00120).
}\medskip
\\
{\small IRMAR (Institut de Recherche Math\'ematique de Rennes,  UMR 6625 du CNRS)}\\
{\small Universit\'e de Rennes 1, campus de Beaulieu, F-35042 Rennes Cedex}\\
{\small laurent.moret-bailly@univ-rennes1.fr}\\
{\small \href{http://perso.univ-rennes1.fr/laurent.moret-bailly/}{http://perso.univ-rennes1.fr/laurent.moret-bailly/}}}
\begin{document}
\selectlanguage{french}
\maketitle
\begin{abstract}
On \'etudie les anneaux (notamment noeth\'eriens) dans lesquels l'ensemble des \'el\'ements non nuls est existentiel positif (r\'eunion finie de projections d'ensembles \og al\-g\'e\-briques\fg). Dans le cas noeth\'erien int\`egre, on montre notamment que cette con\-di\-tion est v\'erifi\'ee pour tout anneau qui n'est pas local hens\'elien, et qu'elle ne l'est pas pour un anneau local hens\'elien excellent qui n'est pas un corps.

Ces r\'esultats apportent au passage une r\'eponse \`a une question de Popescu sur l'approximation forte pour les couples hens\'eliens.
\end{abstract}\selectlanguage{english}
\begin{abstract}
We investigate the rings in which the set of nonzero elements is positive-existential (i.e.\ a finite union of projections of  ``algebraic'' sets). In the case of Noetherian domains, we prove in particular that this condition is satisfied whenever the ring in question is not local Henselian, while it is not satisfied for any excellent local Henselian domain which is not a field.

As a byproduct, we obtain an answer to a question of Popescu on strong approximation for Henselian pairs.
\end{abstract}
\selectlanguage{french}

\noindent\textsl{Classification AMS:} 
03C99, 
13E05, 
13J15, 
13B40, 
11U09. 

\section{Introduction}

\subsection{D\'efinissabilit\'e existentielle.}\label{SsecDefExist} Si $A$ est un anneau (commutatif unitaire) et $r$ un entier naturel, un sous-ensemble $Z$ de $A^r$ est dit [$A$-]\emph{existentiel} (resp. [$A$-]\emph{existentiel positif}) s'il existe une formule $\phi(t_1,\ldots,t_r)$ du langage des anneaux avec symboles de constantes pour les \'el\'ements de $A$,  \`a $r$ variables libres, sans quantificateur universel ($\forall$) (resp. et sans n\'egation) telle que pour tout $\soul{t}\in A^r$, $\phi(\soul{t})$ soit vraie \ssi\ $\soul{t}\in Z$. 

Plus explicitement, un ensemble existentiel (resp. existentiel positif) est r\'eunion d'une famille finie (\'eventuellement vide, voir remarque \ref{RemVide} ci-dessous) d'ensembles de la forme
$$\bigl\{\,\soul{t}\in A^r\mid \exists \soul{x}\in A^n,\;F_1(\soul{t},\soul{x})=\cdots =F_s(\soul{t},\soul{x})=0\wedge G_1(\soul{t},\soul{x})\neq0\wedge\cdots\wedge G_u(\soul{t},\soul{x})\neq0\,\bigr\},
\leqno{(*)}$$
resp.
$$\bigl\{\,\soul{t}\in A^r\mid \exists \soul{x}\in A^n,\;F_1(\soul{t},\soul{x})=\cdots =F_s(\soul{t},\soul{x})=0\,\bigr\}, \leqno{(**)}$$
o\`u les $F_j$ et les $G_j$ sont des  \pols\ en $r+n$ ind\'etermin\'ees \`a \coefs\ dans $A$.

Si $f:A^r\to A^s$ est une application polynomiale (par exemple $A$-lin\'eaire), l'image (resp. l'image r\'eciproque) par $f$ de tout sous-ensemble existentiel de $A^r$ (resp. de $A^s$) est existentielle, et de m\^eme pour les sous-ensembles existentiels positifs. On voit en particulier que la notion de sous-ensemble existentiel (resp. existentiel  positif) a un sens dans tout $A$-module libre de rang fini, ind\'ependamment du choix d'une base. 

\subsubsection{Remarque. }\label{RemVide} Dans la d\'efinition d'un ensemble existentiel positif, la r\'eunion de la \emph{famille vide} d'ensembles $(**)$ donne naissance au sous-ensemble vide de $A^r$, qui est donc existentiel positif. Au niveau des formules, cette op\'eration correspond \`a la \emph{disjonction de la famille vide} de formules, qui est la constante logique \og \textsc{faux}\fg. Il faut donc adjoindre celle-ci au langage des anneaux usuel, ce qui ne semble pas \^etre de pratique courante en th\'eorie des mod\`eles. 

Lorsque $A$ n'est pas nul, la constante \og \textsc{faux}\fg\ est \'equivalente \`a la formule $1=0$, de sorte  que l'on peut s'en dispenser. En revanche, si $A$ est l'anneau nul, la partie vide de $A^r$ n'est pas r\'eunion d'une famille non vide d'ensembles  $(**)$, ceux-ci n'\'etant jamais vides.

Le lecteur pourra, s'il y tient, revenir \`a la d\'efinition traditionnelle $\{+,.,0,1\}$ du langage des anneaux, au prix de modifications mineures (dans certains \'enonc\'es il faut se limiter aux anneaux non nuls).

\subsection{La condition (C).}\label{SsecCondC} Pour $A$ donn\'e, on peut se demander si tout ensemble existentiel dans $A^r$ est existentiel positif. Cette propri\'et\'e \'equivaut manifestement \`a la con\-dition suivante:
\smallskip

\noindent(C)\hfill  \og l'ensemble $A\setminus\{0\}$ est existentiel positif\fg\hfill{ }
\smallskip

\subsubsection{Exemples. }\label{RemCondC}  Il est clair que tout anneau fini v\'erifie (C) (y compris l'anneau nul, vu nos conventions). 

Tout corps $K$ v\'erifie (C) ($t\in K$ est non nul \ssi\ il existe $x\in K$ tel que $tx=1$). 

Il est aussi  bien connu que $\ZZ$ v\'erifie (C): par exemple, si $t\in\ZZ$, alors $t\neq0$ \ssi\ il existe $x,y,w$ dans $\ZZ$ tels que $tw=(1+2x)(1+3y)$. La m\^eme astuce montre d'ailleurs que tout anneau d'entiers alg\'ebriques v\'erifie (C).

L'anneau $A$ des entiers d'un corps local \emph{ne v\'erifie pas} (C): en effet, tout sous-ensemble existentiel positif de $A$ est compact, alors que $A\setminus\{0\}$ ne l'est pas. Plus g\'e\-n\'e\-ra\-le\-ment, un anneau topologique compact infini ne v\'erifie pas (C).

Il est facile de voir qu'un produit $A_1\times A_2$ v\'erifie (C) \ssi\ $A_1$ et $A_2$ v\'erifient (C). En revanche, un \emph{produit infini d'anneaux non nuls} ne v\'erifie jamais (C) (tout ensemble existentiel positif est ferm\'e pour le produit des topologies discr\`etes). 

\subsection{R\'esultats. } Les principaux r\'esultats de ce travail sont les suivants.

\subsubsection{Anneaux noeth\'eriens. }\label{ResultNoeth} Soit $A$ un anneau noeth\'erien. Alors:
\begin{itemize}
\item si $A$ est int\`egre et n'est pas local hens\'elien, il v\'erifie \textup{(C)} (\ref{ThNonHens});
\item si $A$ est un localis\'e d'un anneau de Jacobson noeth\'erien, il v\'erifie (C) (\ref{CorJac});
\item si $A$ est local hens\'elien excellent de dimension $>0$, il ne v\'erifie pas \textup{(C)} (\ref{ThHens}).
\end{itemize}
\subsubsection{Propri\'et\'es d'approximation. }Ces propri\'et\'es (notamment les r\'esultats de \cite{PfiPop75,DeLi80,Pop86,Swan95, Spiv99}) jouent un r\^ole essentiel dans la preuve de \ref{ThHens}; \emph{a contrario}, on d\'eduit de \ref{ThNonHens} que si un couple hens\'elien $(A,\fr)$ v\'erifie la \og propri\'et\'e d'approximation forte\fg, alors $A$ est n\'ecessairement semi-local hens\'elien (sauf cas triviaux, comme $\fr=0$). Ce r\'esultat (corollaire \ref{CorBiz}) pr\'ecise la r\'eponse n\'egative donn\'ee par Spivakovsky \cite{Spiv94} \`a une question de D. Popescu; voir la remarque \ref{RemCorBiz}.
\subsubsection{Anneaux de fonctions holomorphes. }\label{ResultHolo} Soit $X$ un espace analytique de Stein (par exemple un ferm\'e analytique de $\CC^n$ ou d'un polydisque ouvert), irr\'eductible et r\'eduit. Alors l'anneau des fonctions holomorphes sur $X$ v\'erifie (C) (th\'eor\`eme \ref{ThStein}).

\subsection{Remarques sur les m\'ethodes. } Les deux \'el\'ements essentiels dans la d\'emonstration du th\'eor\`eme \ref{ThNonHens} sont le lemme \ref{LemDeuxId}, qui suppose l'existence de deux id\'eaux premiers ayant certaines propri\'et\'es, et le lemme \ref{LemmeDoublement} qui permet de remplacer l'anneau \'etudi\'e (muni d'un id\'eal premier $\fp$) par un autre ayant deux id\'eaux premiers au-dessus de $\fp$. 

De ces deux lemmes, le premier g\'en\'eralise \`a peu de chose pr\`es l'argument utilis\'e plus haut pour $\ZZ$ (\ref{RemCondC}, o\`u les id\'eaux en question sont $2\ZZ$ et $3\ZZ$), dont des variantes ont \'et\'e maintes fois utilis\'ees dans la litt\'erature (voir par exemple \cite{DavMatRob74} ou, pour une g\'en\'eralisation r\'ecente, \cite{Dem07}, Proposition 2.3); l'id\'ee du second remonte au moins \`a \cite{Shlap94}, Theorem 4.2.

Quant au th\'eor\`eme \ref{ThHens}, il repose comme on l'a dit sur les propri\'et\'es d'approximation, le lien \'etant la proposition \ref{PropFerme}. Le lien entre la condition (C) et l'approximation \'etait d\'ej\`a connu de T. Pheidas. 

On voit donc que ce travail ne contient pas d'id\'ee essentiellement nouvelle, l'auteur s'\'etant seulement efforc\'e de syst\'ematiser les m\'ethodes existantes.

\section{G\'en\'eralit\'es et rappels.}\label{SecRapp}

\subsection{D\'efinissabilit\'e.} 

On laisse au lecteur la d\'emonstration du lemme \ref{QuotientDioph} ci-dessous.

\begin{sublem}\label{QuotientDioph} Soient $A$ un anneau, $I$ un id\'eal de $A$, $\pi:A\to A/I$ l'homomorphisme canonique.
\begin{romlist}
\item\label{QuotientDioph1}  Si $I$ est de type fini, c'est un sous-ensemble existentiel positif de $A$.
\item\label{QuotientDioph2}  On suppose que $I$ est existentiel positif  dans $A$. Si $D$ est un sous-ensemble existentiel positif de $A/I$, alors $\pi^{-1}(D)$ est existentiel positif dans $A$. 

En particulier, si $A/I$ v\'erifie \textup{(C)}, alors $A\setminus I$ est existentiel positif dans $A$.\qed
\end{romlist}
\end{sublem}

\begin{sublem}\label{LemWeil} Soient $A$ un anneau, et $B$ une $A$-alg\`ebre qui est  un $A$-module libre de rang fini $d$.
\begin{romlist}
\item\label{LemWeil1}  Soit $Z$ un sous-ensemble $B$-existentiel positif de $B$. Alors $Z$ est $A$-existentiel positif (dans $B$ vu comme $A$-module).
\item\label{LemWeil2}  Si $d>0$ et si $B$ v\'erifie \textup{(C)}, alors $A$ v\'erifie \textup{(C)}.
\end{romlist}
\end{sublem}
\dem on d\'eduit \ref{LemWeil1} de la remarque suivante: si $F:B^r\to B$ est d\'efinie par un \pol\ \`a \coefs\ dans $B$, alors $F$  est d\'efinie par $d$ \pols\ (en $rd$ ind\'etermin\'ees) \`a \coefs\ dans $A$, une fois $B$ identifi\'e \`a $A^d$ par le choix d'une base. On en d\'eduit \ref{LemWeil2} en prenant $Z=B\setminus\{0\}$ et en remarquant que $A\setminus\{0\}$ est image r\'eciproque de $Z$ par l'application $A$-lin\'eaire canonique de $A$ dans $B$, qui est injective si $d>0$.\qed

\subsection{Anneaux locaux hens\'eliens.}\label{SsecLocHens}

On rappelle \cite{EGA4_IV,RayHens} qu'un anneau local $A$, de corps r\'esiduel $k$, est dit \emph{hens\'elien} si toute $A$-alg\`ebre finie (\emph{i.e.} de type fini comme $A$-module) est un produit d'anneaux locaux. Une d\'efinition \'equivalente est la suivante: pour tout \pol\ unitaire $P\in A[X]$ et tout \'el\'ement $\ol{x}$ de $k$ qui est racine simple de (l'image de) $P$, il existe un unique $x\in A$ relevant $\ol{x}$ et annulant $P$. \smallskip

Pour utiliser efficacement la premi\`ere d\'efinition, nous aurons besoin du lemme suivant:

\begin{sublem}\label{LemProdLoc} Soit $A$ un anneau semi-local noeth\'erien. Les conditions suivantes sont \'equivalentes:
\begin{romlist}
\item\label{LemProdLoc1} $A$ est isomorphe \`a un produit d'anneaux locaux;
\item\label{LemProdLoc2} tout id\'eal premier de $A$ est contenu dans un \emph{unique} id\'eal maximal (autrement dit, tout quotient int\`egre de $A$ est local);
\item\label{LemProdLoc3} tout id\'eal premier minimal de $A$ est contenu dans un unique id\'eal maximal.
\end{romlist}
\end{sublem}
\dem les implications \ref{LemProdLoc1}\,$\Rightarrow$\ref{LemProdLoc2}\,$\Leftrightarrow$\ref{LemProdLoc3} sont faciles et laiss\'ees au lecteur. Supposons \ref{LemProdLoc3}, et montrons \ref{LemProdLoc1}. Posons $X=\mathrm{Spec}\,A$, et soient $x_1,\ldots,x_r$ les points ferm\'es de $X$ (correspondant aux id\'eaux maximaux $\fm_1,\ldots,\fm_r$ de $A$). Pour chaque $i\in\{1,\ldots,r\}$, notons $X_i$ le ferm\'e de $X$ r\'eunion des composantes \irrs\ contenant $x_i$. L'hypoth\`ese \ref{LemProdLoc3} signifie que $X$ est r\'eunion \emph{disjointe} des $X_i$, qui sont donc ouverts dans $X$. Si l'on munit chaque $X_i$ de sa structure de sous-sch\'ema ouvert de $X$, alors $X_i$ est local de point ferm\'e $x_i$ et s'identifie donc \`a $\mathrm{Spec}\,A_{\fm_i}$. Le fait que $X$ soit r\'eunion disjointe des $X_i$ entra\^{\i}ne donc que $A\cong\prod_{i=1}^r A_{\fm_i}$, d'o\`u la conclusion.\qed

\begin{subcor}\label{CorProdLoc} Soit $A$ un anneau local noeth\'erien. Les conditions suivantes sont \'equi\-va\-lentes:
\begin{romlist}
\item\label{CorProdLoc1} $A$ est hens\'elien;
\item\label{CorProdLoc2} toute $A$-alg\`ebre finie \emph{int\`egre} est un  anneau local.\qed
\end{romlist}
\end{subcor}

\subsection{Anneaux excellents.} On renvoie \`a \cite{EGA4_II} ou \cite{Matsu89} pour la d\'efinition des anneaux excellents. Rappelons seulement quelques propri\'et\'es:
\begin{romlist}
\item tout anneau excellent est noeth\'erien;
\item tout anneau local noeth\'erien complet est excellent;
\item si $A$ est excellent, il en est de m\^eme de toute $A$-alg\`ebre de type fini et de tout anneau de fractions de $A$;
\item si $A$ est local excellent, il en est de m\^eme de son hens\'elis\'e;
\item un anneau de valuation discr\`ete $A$ est excellent \ssi\ le corps des fractions de son compl\'et\'e est une extension s\'eparable du corps des fractions de $A$.
\end{romlist}

\subsection{Propri\'et\'es d'approximation.}\label{PropAppr}

Soient $A$ un anneau noeth\'erien et $\fr$ un id\'eal de $A$. Pour tout $q\in\NN$, on pose $A(q):=A/\fr^{q}$, et l'on note $\wh{A}=\varprojlim_{q\in\NN} A(q)$ le s\'epar\'e compl\'et\'e $\fr$-adique de $A$. 

Pour toute famille finie $S=(F_j)_{j=1,\ldots,s}$ de \pols\ en $r$ ind\'etermin\'ees $\uple{r}$ \`a \coefs\ dans $A$, et pour toute $A$-alg\`ebre $B$, on posera 
$$\mathrm{sol}(S,B):=\{\soul{x}\in B^r\mid F_j(\soul{x})=0 \text{ pour tout }j\}.$$
Il est clair que $\mathrm{sol}(S,B)$ est fonctoriel en $B$; il s'identifie \`a l'ensemble des morphismes de $A$-alg\`ebres de $A[\soul{X}]/(\uple[F]{s})$ dans $B$. 

On consid\'erera la condition suivante (\og principe de Hasse infinit\'esimal\fg):
\smallskip

\noindent PHI($A,\fr$): \og pour toute famille $S$ comme ci-dessus, si $\mathrm{sol}(S,A(q))\neq\emptyset$ pour tout $q\in\NN$ alors $\mathrm{sol}(S,{A})\neq\emptyset$\fg.

\subsubsection{Remarques. }\label{RemPHI} (1) Il est facile de voir que PHI($A,\fr$) implique les propri\'et\'es suivantes:
\begin{romlist}
\item\label{RemPHI1} $\fr$ est contenu dans le radical de $A$;
\item\label{RemPHI2} (cons\'equence de \ref{RemPHI1}) $A$ est $\fr$-adiquement s\'epar\'e, \emph{i.e.} $\bigcap_{q\in\NN}\fr^q =\{0\}$;
\item\label{RemPHI3} pour tout id\'eal $J$ de $A$, PHI($A/J,\,(\fr+J)/J$) est v\'erifi\'ee;
\item\label{RemPHI4} $(A,\fr)$ v\'erifie la \og propri\'et\'e d'approximation\fg: pour $S$ comme ci-dessus, $\mathrm{sol}(S,A)$ est dense dans $\mathrm{sol}(S,\wh{A})$ pour la topologie $\fr$-adique;
\item\label{RemPHI5} (cons\'equence de \ref{RemPHI4}) $(A,\fr)$ est un \emph{couple hens\'elien}: si $B$ est une $A$-alg\`ebre \'etale telle que $A/\fr\flis B/\fr B$, alors il existe un $A$-homomorphisme de $B$ dans $A$.
\end{romlist}
(2) D'autre part, PHI($A,\fr$) est vraie d\`es que $(A,\fr)$ v\'erifie la \og propri\'et\'e d'approximation forte\fg, que l'on peut formuler ainsi: pour $S$ comme ci-dessus et tout $q\in\NN$, il existe $q'\geq q$ tel que $\mathrm{sol}(S,A)$ et $\mathrm{sol}(S,A(q'))$  aient m\^eme image dans $\mathrm{sol}(S,A(q))$.
\smallskip

\noindent(3) Lorsque $A$ est local d'id\'eal maximal $\fm$, la propri\'et\'e d'approximation pour $(A,\fm$) est en fait \emph{\'equivalente} \`a la propri\'et\'e d'approximation forte (et donc aussi \`a PHI($A,\fm$)), d'apr\`es \cite{PfiPop75} (red\'emontr\'e dans  \cite{DeLi80} par des m\'ethodes de th\'eorie des mod\`eles). 
\smallskip

\noindent(4) Enfin, lorsque $A$ est local hens\'elien \emph{excellent} d'id\'eal maximal $\fm$, la propri\'et\'e d'approximation est v\'erifi\'ee, d'apr\`es \cite{Pop86} (autres r\'ef\'erences:  \cite{Swan95}, \cite{Spiv99}). Il en est donc de m\^eme de PHI($A,\fm$), et \emph{a fortiori} de PHI($A,\fr$) pour tout id\'eal strict $\fr$. Retenons donc pour la suite:

\begin{thm}\label{ThAppr}  Soient $A$ un anneau \emph{local hens\'elien excellent}, et $\fr$ un id\'eal strict de $A$. Alors  la condition \textup{PHI($A,\fr$)} est v\'erifi\'ee.\qed
\end{thm}

 \section{Anneaux noeth\'eriens int\`egres}

Dans ce paragraphe, nous allons \'etablir le th\'eor\`eme suivant:
\begin{thm}\label{ThNonHens} Soit $A$ un anneau int\`egre noeth\'erien. Si $A$ n'est pas local hens\'elien, alors $A$ v\'erifie \textup{(C)}.
\end{thm}

\begin{lem}\label{LemDeuxId} Soient $A$ un anneau int\`egre, $\fp_1$ et $\fp_2$ deux id\'eaux premiers de $A$. On suppose que:
\begin{romlist}
\item\label{LemDeuxId1} $A_{{\fp}_1}$ est noeth\'erien;
\item\label{LemDeuxId2} $\fp_1\cap\fp_2$ ne contient aucun id\'eal premier non nul de $A$.
\end{romlist}
Alors, pour tout $t\in A$, on a l'\'equivalence:
$$t\neq0\;\Longleftrightarrow \;\exists(w,x_1,x_2)\in A^3,\:
\begin{cases}
t\,w=x_1\,x_2\\
x_1\not\in\fp_1\\
x_2\not\in\fp_2.
\end{cases}
$$
En outre, si $\fp_1$ et $\fp_2$ sont de type fini et si $A/\fp_1$ et $A/\fp_2$ v\'erifient \textup{(C)}, alors $A$ v\'erifie \textup{(C)}.
\end{lem}
\dem
 la derni\`ere assertion r\'esulte de la premi\`ere en vertu du lemme \ref{QuotientDioph} appliqu\'e aux id\'eaux $\fp_1$ et $\fp_2$. 

Montrons la premi\`ere assertion. L'implication $\Leftarrow$ est triviale puisque $A$ est int\`egre. R\'eciproquement, soit $t\neq0$ dans $A$. Comme $A_{\fp_1}$ est noeth\'erien, les id\'eaux premiers minimaux de $A_{\fp_1}/(t)$ sont en nombre fini; ils correspondent \`a  un nombre fini d'id\'eaux premiers de $A$ contenant $t$ et contenus dans $\fp_1$, que nous noterons $\fq_1,\ldots,\fq_r$. L'hypoth\`ese \ref{LemDeuxId2} implique qu'aucun des $\fq_j$ n'est contenu dans $\fp_2$, ce qui implique que $\bigcap_{j}\fq_j\not\subset\fp_2$. Soit donc $y\in\bigl(\bigcap_{j}\fq_j\bigr)\setminus\fp_2$. L'image de $y$ dans $A_{\fp_1}/(t)$ appartient \`a tous les id\'eaux premiers minimaux donc est nilpotente: il existe $n\in\NN$ et $v\in A_{\fp_1}$ tels que $y^n=tv$. Vu la d\'efinition du localis\'e $A_{\fp_1}$, il existe donc $s\in A\setminus\fp_1$ et $w\in A$ tels que $sy^n=tw$, de sorte que $t$ v\'erifie la propri\'et\'e voulue (avec $x_1=s$ et $x_2=y^n$: noter que $y^n\not\in\fp_2$ puisque $y\not\in\fp_2$).\qed
\subsubsection{Remarques. }\label{RemLemDeuxId}{\ }\smallskip

\noindent(1) La condition \ref{LemDeuxId1} de l'\'enonc\'e peut \^etre affaiblie en \og l'espace topologique $\Spec(A_{\fp_1})$ est noeth\'erien\fg; plus pr\'ecis\'ement on n'utilise que le fait que pour tout $t\in A_{\fp_1}$, l'anneau $A_{\fp_1}/(t)$ n'a qu'un nombre fini d'id\'eaux premiers minimaux.\smallskip

\noindent(2) La condition \ref{LemDeuxId2} est notamment v\'erifi\'ee si l'on a  $\dim A_{\fp_1}=1$ et $\fp_1\not\subset\fp_2$ (ou l'inverse).\medskip

Voici un corollaire imm\'ediat de \ref{LemDeuxId}:
\begin{subcor}\label{CorA[X]} Soit $R$ un anneau int\`egre v\'erifiant  \textup{(C)}.
Alors $R[X]$ v\'erifie  \textup{(C)}.
\end{subcor}
\dem il suffit d'appliquer la derni\`ere assertion de \ref{LemDeuxId} avec $A=R[X]$, $\fp_1=(X)$ et $\fp_2=(X-1)$.\qed

\begin{lem}\label{LemmeDoublement} Soit $A$ un anneau int\`egre noeth\'erien, de corps des fractions $K$. Soit $\fp$ un id\'eal premier non nul de $A$; on suppose que $\fp$ n'est pas le plus grand id\'eal premier de $A$ (autrement dit, $A$ n'est pas un anneau local d'id\'eal maximal $\fp$).

Alors il existe un \pol\ $F=X^2+aX+b\in A[X]$ ayant les propri\'et\'es suivantes:
\begin{romlist}
\item\label{LemmeDoublement1} $a\not\in\fp$;
\item\label{LemmeDoublement2} $b\in\fp$;
\item\label{LemmeDoublement3} $F$ est irr\'eductible dans $K[X]$.
\end{romlist}
\end{lem}
\dem l'hypoth\`ese sur $\fp$ entra\^{\i}ne qu'il existe un \'el\'ement non inversible $\alpha$ de $A\setminus\{0\}$ qui n'appartient pas \`a $\fp$. Soit $\fq$ un id\'eal premier  minimal parmi ceux contenant $\alpha$: alors $\fq\not\subset\fp$ (puisque $\alpha\in\fq$), et l'anneau $A_\fq$ est local noeth\'erien de dimension $1$ (vu la minimalit\'e de $\fq$). Comme $\fp$ n'est pas nul cela entra\^{\i}ne que $\fp\not\subset\fq$. Fixons $\beta\in\fp\setminus\fq$. Nous avons ainsi trouv\'e un id\'eal premier $\fq$ et deux \'el\'ements $\alpha$ et $\beta$ de $A$ v\'erifiant:\begin{center}
$\alpha\in\fq\setminus\fp$; $\beta\in\fp\setminus\fq$; $\dim A_\fq=1$.
\end{center}
Il r\'esulte du th\'eor\`eme de Krull-Akizuki (\cite{BourAC7}, \S2, n\up{\accent'27} 5, prop. 5 et cor. 2) que la cl\^oture int\'egrale de $A_\fq$ est un anneau de Dedekind; si on le localise en l'un de ses id\'eaux maximaux, on obtient un anneau de valuation discr\`ete de corps des fractions $K$ qui domine $A_\fq$. Notons $\nu:K\to\ZZ\cup\{+\infty\}$ la valuation normalis\'ee correspondante. Comme le mono\"{\i}de $\nu(A\setminus\{0\})$ engendre le groupe $\ZZ$, il existe $\gamma\in A$ tel que $\nu(\gamma)$ soit impair et positif. 

Posons $b:=\beta\gamma$. Alors $b\in\fp$ et $\nu(b)$ est impair. Pour $n\in\NN$ convenable, l'\'el\'ement $a:=\alpha^n$ v\'erifie $\nu(a)>\nu(b)/2$ (et \'evidemment $a\not\in\fp$). 

Il reste \`a montrer  que $F=X^2+aX+b$ n'a pas de racine dans $K$. Si $z$ \'etait une telle racine, on aurait $z(z+a)=-b$. Si $\nu(z)<\nu(a)$ ceci entra\^{\i}ne $\nu(z(z+a))=2\nu(z)=\nu(b)$, impossible puisque $\nu(b)$ est impair; sinon, on a 
$\nu(z)\geq\nu(a)$ donc  $\nu(z(z+a))\geq2\nu(a)>\nu(b)$ vu le choix de $a$: nouvelle contradiction. (Les amateurs de polygones de Newton se contenteront de remarquer que celui de $F$ a pour seule pente $-\nu(b)/2$, qui n'est pas un entier).\qed

\subsubsection{Remarque. }\label{RemLemmeDoublement} On voit notamment que, sous les hypoth\`eses de \ref{LemmeDoublement}, l'anneau local  $A_\fp$ n'est pas hens\'elien, puisque $F$ est \irr\ mais a deux racines simples dans le corps r\'esiduel de $A_\fp$.

\begin{prop}\label{PropRec} Soit $A$ un anneau int\`egre noeth\'erien. On suppose qu'il existe un id\'eal premier $\fp$ de $A$ tel que:
\begin{romlist}
\item\label{PropRec1} $\fp$ n'est pas le plus grand id\'eal premier de $A$;
\item\label{PropRec2} $A/\fp$ v\'erifie \textup{(C)}.
\end{romlist}
Alors  $A$ v\'erifie \textup{(C)}.
\end{prop}
\dem on proc\`ede par r\'ecurrence sur $h:=\dim A_\fp$. Si $h=0$, alors $\fp$ est nul et l'assertion est triviale.

Supposons $h>0$. Il existe alors un \pol\ $F=X^2+aX+b$ comme dans le lemme \ref{LemmeDoublement}. Soit $B$ la $A$-alg\`ebre $A[X]/(F)$. Alors $B$ est int\`egre, et est un $A$-module libre de rang 2; il suffit donc (lemme \ref{LemWeil}) de montrer que $B$ v\'erifie (C). L'anneau $B/\fp B$ est isomorphe \`a $(A/\fp)[X]/(X(X+\ol{a}))$ o\`u $\ol{a}\neq0$ est la classe de $a$ modulo $\fp$. Donc $B$ a deux id\'eaux premiers distincts au-dessus de $\fp$, \`a savoir $\fp_1=\fp B+ xB$ et $\fp_2=\fp B+(x+a)B$ o\`u $x$ est la classe de $X$.   Comme $B$ est fini libre sur $A$, on a $\dim B_{\fp_1}=\dim B_{\fp_2}=h$. De plus $B/\fp_1\cong A/\fp$ (l'isomorphisme envoie la classe de $x$ sur $0$) et de m\^eme $B/\fp_2\cong A/\fp$ (l'isomorphisme envoie la classe de $x$ sur $-\ol{a}$). En particulier les anneaux $B/\fp_1$ et $B/\fp_2$ v\'erifient (C). Distinguons deux cas:
\smallskip

\noindent(1) $\fp_1\cap\fp_2$ ne contient aucun id\'eal premier non nul de $B$ (condition v\'erifi\'ee en particulier si $h=1$): alors on d\'eduit du lemme \ref{LemDeuxId} que $B$ v\'erifie (C). \smallskip

\noindent(2) il existe un id\'eal premier $\fq\subset\fp_1\cap\fp_2$ non nul: en choisissant $\fq$ minimal on peut supposer que $\dim B_\fq=1$. Pour $i\in\{1,2\}$ posons $\ol{\fp}_i:=\fp_i/\fq\subset \ol{B}:=B/\fq$. Alors les localis\'es  $\ol{B}_{\ol{\fp}_i}$ sont de dimension $\leq h-1$; en outre $\ol{\fp}_1\not\subset \ol{\fp}_2$, donc $\ol{\fp}_2$ n'est pas le plus grand id\'eal premier de $\ol{B}$ et l'hypoth\`ese de r\'ecurrence s'applique \`a $\ol{B}$. Ainsi $\ol{B}$  v\'erifie (C), et donc $B$ aussi d'apr\`es le cas $h=1$ d\'ej\`a \'etabli.\qed

\subsection{D\'emonstration du th\'eor\`eme \ref{ThNonHens}. } Soit $A$ comme dans l'\'enonc\'e. Si $A$ n'est pas local, il suffit d'appliquer \ref{PropRec}
en prenant pour $\fp$ n'importe quel id\'eal maximal de $A$: il est clair que la condition \ref{PropRec2} de la proposition est v\'erifi\'ee puisque $A/\fp$ est un corps.

Supposons  $A$ local mais non hens\'elien. Il existe alors une $A$-alg\`ebre finie int\`egre $B_0$ qui n'est pas locale, d'apr\`es \ref{CorProdLoc}.  
Soit $(\uple[\xi]{s})$ une famille g\'en\'eratrice finie du $A$-module $B_0$, et pour chaque $i\in\{1,\ldots,s\}$ soit $F_i\in A[X_i]$ un \pol\ unitaire, en une ind\'etermin\'ee $X_i$, annulant $\xi_i$. Alors $B_0$ est quotient de $B:=A[\uple{s}]/(\uple[F]{s})$ qui est libre de rang fini comme $A$-module. 

Nous avons donc trouv\'e une $A$-alg\`ebre finie libre $B$ et un id\'eal premier $\fp$ de $B$ tel que $B_0=B/\fp$ ne soit pas local. Soit  $\fq$ un id\'eal premier \emph{minimal} de $B$ contenu dans $\fp$: alors $B/\fq$ n'est pas local donc v\'erifie (C). Par suite $B\setminus\fq$ est existentiel positif dans $B$ (\ref{QuotientDioph}). Soit $j:A\to B$ le morphisme structural. Comme les \'el\'ements de $\fq$ sont diviseurs de z\'ero dans $B$, que $A$ est int\`egre et que $B$ est $A$-libre on a $j^{-1}(\fq)=\{0\}$. Donc $A\setminus\{0\}=j^{-1}(B\setminus\fq)$ est existentiel positif dans $A$, ce qui ach\`eve la d\'emonstration.\qed 

\subsection{Anneaux de fractions. } Une cons\'equence du th\'eor\`eme \ref{ThNonHens} est que pour les anneaux noeth\'eriens int\`egres, la propri\'et\'e (C) est \og stable par anneaux de fractions\fg. Plus pr\'ecis\'ement:

\begin{subcor}\label{CorLocalise} Soit $A$ un anneau int\`egre noeth\'erien, et soit $S$ une partie multiplicative de $A$. On suppose v\'erifi\'ee l'une des conditions suivantes:
\begin{romlist}
\item\label{CorLocalise1} $S\not\subset A^\times$ \emph{(de sorte que $S^{-1}A\neq A$)};
\item\label{CorLocalise2} $A$ v\'erifie \textup{(C)}.
\end{romlist} 
Alors $S^{-1}A$ v\'erifie \textup{(C)}.
\end{subcor}

\dem Si $0\in S$, alors $S^{-1}A$ est nul et tout est trivial. On  suppose donc que $0\not\in S$,  de sorte que $S^{-1}A$ est noeth\'erien et int\`egre. S'il n'est pas local, on conclut par \ref{ThNonHens}. Supposons-le  local:  il est donc de la forme $A_\fp$ o\`u $\fp$ est un id\'eal premier de $A$. Si $\fp$ est nul, alors $S^{-1}A$ est un corps; si $A$ est local d'id\'eal maximal $\fp$, alors $S^{-1}A=A$, ce qui est  exclu dans le cas \ref{CorLocalise1} et implique le r\'esultat dans le cas \ref{CorLocalise2}. Sinon, nous sommes dans la situation du lemme \ref{LemmeDoublement}, de sorte que  $A_\fp$ n'est pas hens\'elien d'apr\`es \ref{RemLemmeDoublement}, donc v\'erifie (C).\qed

\section{Le cas hens\'elien}\label{SecHens}
\begin{thm}\label{ThHens} Soient $A$ un anneau noeth\'erien, $\fr$ un id\'eal non contenu dans le nilradical de $A$.  Si  la propri\'et\'e \textup{PHI($A,\fr$)} de \textup{\ref{PropAppr}} est v\'erifi\'ee
(et notamment si $A$ est local hens\'elien excellent, d'apr\`es \textup{\ref{ThAppr}}), alors $A$ ne v\'erifie pas \textup{(C)}.
\end{thm}
 \dem l'hypoth\`ese sur $\fr$  \'equivaut \`a dire que la topologie $\fr$-adique de $A$ n'est pas discr\`ete. Donc $A\setminus\{0\}$ n'est pas ferm\'e dans $A$ pour cette topologie, de sorte que le th\'eor\`eme r\'esulte de la proposition  \ref{PropFerme} qui suit.\qed
 
 \begin{subprop}\label{PropFerme} Soit $A$ un anneau noeth\'erien, $\fr$ un id\'eal de $A$.   Les conditions suivantes sont \'equivalentes:
\begin{romlist}
\item\label{PropFerme1} \textup{PHI($A,\fr$)} est v\'erifi\'ee;
\item\label{PropFerme2} pour tout $n\in\NN$, toute partie existentielle positive de $A^n$ est ferm\'ee pour la topologie $\fr$-adique.
\end{romlist}
\end{subprop} 

\dem \ref{PropFerme1}\,$\Rightarrow$\,\ref{PropFerme2}: supposons {PHI($A,\fr$)} v\'erifi\'ee, et soit $Z\subset A^n$ un sous-ensemble existentiel positif. Alors $Z$ est r\'eunion finie d'ensembles de la forme
$$\bigl\{\,\soul{t}\in A^n \mid 
(\exists\,\soul{x}) \;F_{1}(\soul{t},\soul{x})=\ldots=F_{s}(\soul{t},\soul{x})=0\,\bigr\}$$
o\`u chaque  $F_{j}$ est un \pol\ \`a \coefs\ dans $A$ en $n+r$ variables (o\`u $\soul{x}=(\uple[x]{r})$). Pour voir que $Z$ est ferm\'e, on peut donc supposer qu'il est de la forme ci-dessus. 

Soit $\soul{t}\in A^n$ adh\'erent \`a $Z$. Pour tout $q\in\NN$, il existe donc $(\soul{t}(q),\soul{x}(q))\in A^{n+r}$ tels que $\soul{t}-\soul{t}(q)\in \fr^{q}A^{n+r}$ et que $F_j(\soul{t}(q),\soul{x}(q))=0$ pour tout $j$. Comme les $F_j$ sont \`a \coefs\ dans $A$, on a donc $F_j(\soul{t},\soul{x}(q))\in\fr^{q}$. Autrement dit, le syst\`eme d'\'equations 
$$F_j(\soul{t},X_1,\ldots,X_r)=0\qquad(j=1,\ldots,s)$$ 
en les inconnues $X_i$ admet une solution modulo $\fr^{q}$ pour tout $q$. L'hypoth\`ese implique donc qu'il a une solution dans $A^r$, donc que $\soul{t}\in Z$, cqfd.
\smallskip

\ref{PropFerme2}\,$\Rightarrow$\,\ref{PropFerme1}: supposons \ref{PropFerme2}, et reprenons les notations de \ref{PropAppr}. On a donc une famille finie $S=(F_j)_{j=1,\ldots,s}$ de \pols\ en $r$ ind\'etermin\'ees \`a \coefs\ dans $A$. Supposons que $\mathrm{sol}(S,A(q))\neq\emptyset$ pour tout $q\in\NN$, et consid\'erons l'ensemble $Z\subset A^s$ image de l'application $(\uple[F]{s}):A^r\to A^s$. Alors $Z$ est existentiel positif, et l'hypoth\`ese sur $S$ signifie que $0\in A^s$ est adh\'erent \`a $Z$. Comme $Z$ est ferm\'e d'apr\`es \ref{PropFerme2}, on conclut que $0\in Z$, donc que $\mathrm{sol}(S,A)\neq\emptyset$.\qed\medskip

\subsection{Application \`a une question de D. Popescu. }\label{QuestPop}
En combinant \ref{ThHens} et \ref{ThNonHens}, on obtient:

\begin{subcor}\label{CorBiz} Soient $A$ un anneau int\`egre noeth\'erien, $\fr$ un id\'eal non nul de $A$. Si \textup{PHI($A,\fr$)} est v\'erifi\'ee, alors $A$ est local hens\'elien.\qed
\end{subcor}

\subsubsection{Remarque. }\label{RemCorBiz} En utilisant la remarque \ref{RemPHI}, on en d\'eduit plus g\'en\'eralement que si $A$ est un anneau noeth\'erien, $\fr$ un id\'eal de $A$ non contenu dans un id\'eal premier minimal, et si \textup{PHI($A,\fr$)} est v\'erifi\'ee, alors $A$ est semi-local hens\'elien (donc produit d'anneaux locaux hens\'eliens).

Popescu demandait dans \cite{Pop86} si tout couple hens\'elien $(A,\fr)$ (\emph{cf.} \ref{RemPHI}\,(1)\ref{RemPHI5}) tel que l'homomorphisme de compl\'etion $A\to\wh{A}$ soit r\'egulier v\'erifie l'approximation forte, et Spivakovsky \cite{Spiv94}  a montr\'e par un exemple que la r\'eponse est en g\'en\'eral n\'egative. Le r\'esultat ci-dessus montre que la r\'eponse est \emph{toujours} n\'egative, sauf dans le cas semi-local (o\`u elle est affirmative, au moins lorsque $\fr$ est le radical; j'ignore ce qui se passe  par exemple si $A$ est local hens\'elien excellent et si $\fr$ est diff\'erent de l'id\'eal maximal).
\medskip

On a une r\'eciproque partielle \`a \ref{ThHens}:
\begin{thm}\label{ThHensRec} Soit $A$ un anneau local noeth\'erien hens\'elien, \emph{int\`egre et de dimension $1$}, d'id\'eal maximal $\fm$. Si \textup{PHI($A,\fm$)} n'est pas v\'erifi\'ee, alors $A$ v\'erifie \textup{(C)}. 
\end{thm}

\dem par hypoth\`ese il existe un syst\`eme $S=(F_j)_{j=1,\ldots,s}$ comme dans \ref{PropAppr} qui a un z\'ero modulo $\fm^q$ pour tout $q\in\NN$ mais n'a pas de z\'ero dans $A^r$. 

Montrons qu'il existe un tel syst\`eme r\'eduit \`a un seul \pol: en effet, comme $\dim A>0$, le corps des fractions $K$ de $A$ n'est pas alg\'ebriquement clos (il admet une valuation discr\`ete non triviale, cf. la preuve de \ref{LemmeDoublement}) et il existe donc un \pol\ $P\in K[\uple{s}]$ ayant $(0,\ldots,0)$ pour seul z\'ero dans $K^s$. On peut naturellement prendre $P$ dans $A[\soul{X}]$; le \pol\ compos\'e $F:=P(\uple[F]{s})$ a alors un z\'ero modulo $\fm^q$ pour tout $q\in\NN$ mais n'a pas de z\'ero dans $A^r$. 

Montrons alors que pour $x\in A$ quelconque on a l'\'equivalence:
$$x\neq0\quad\Longleftrightarrow\quad(\exists t)(\exists w)\; xw=F(t).$$
L'implication $\Leftarrow$ r\'esulte du fait que $F$ n'a pas de z\'ero dans $A^r$; r\'eci\-pro\-que\-ment, si $x\neq0$, l'anneau $A/xA$ est artinien donc il existe $q\in\NN$ tel que $\fm^q\subset xA$; vu le choix de $F$, il existe $t\in A$ tel que $F(t)\in\fm^q$, donc $F(t)\in xA$, cqfd.\qed

\subsection{Un exemple.}

Nous donnons ici un exemple, tir\'e de \cite{BoLuRa90}, d'anneau de valuation discr\`ete hens\'elien qui ne v\'erifie pas (PHI). Soit $k$ un corps de caract\'eristique $p>0$. Choisissons un \'el\'ement $\xi$ de $k[[T]]$ qui est transcendant sur $k(T)$; posons $U=\xi^p$, et $A_0:=k[[T]]\cap k(T,U)$. Enfin soit $A$ le hens\'elis\'e de $A_0$. On v\'erifie facilement que $A_0$ et $A$ sont des anneaux de valuation discr\`ete, de compl\'et\'e $k[[T]]$. De plus, par construction, $\xi^p\in A_0$ mais $\xi\not\in A_0$ et donc $\xi\not\in A$ puisque le corps des fractions de $A$ est s\'eparable sur celui de $A_0$. L'\'equation $X^p=U$ n'a donc pas de solution dans $A$, mais en a une dans $A/(T^q)$ pour tout $q>0$ (\`a savoir la classe d'un \pol\ en $T$ congru \`a $\xi$ modulo $T^q$). 

En particulier $A$ v\'erifie (C): explicitement,  pour tout $f\in A$, on a $f\neq0$ \ssi\ il existe $g$ et $h$ dans $A$ tels que $fg=h^p-U$.

\section{Anneaux noeth\'eriens non int\`egres}\label{SecNonInt}

\subsection{Notations. }\label{IntroNonInt} La condition (C), \`a la connaissance de l'auteur, n'a \'et\'e utilis\'ee dans la litt\'erature que dans le cas int\`egre; pour les anneaux plus g\'en\'eraux, nous nous contenterons donc de quelques remarques \'el\'ementaires. Notons d'abord que pour un anneau non int\`egre il est naturel de consid\'erer, outre $A\setminus\{0\}$, l'ensemble  $A^{\mathrm{reg}}$ des \'el\'ements \emph{r\'eguliers} (c'est-\`a-dire \og non diviseurs de z\'ero\fg) de $A$, qui est contenu dans $A\setminus\{0\}$ si $A$ n'est pas nul et lui est \'egal si $A$ est int\`egre. Ceci conduit \`a envisager, pour un anneau $A$ donn\'e, la condition
\smallskip

\noindent$\mathrm{(C^{reg})}$\hfill  \og l'ensemble $A^{\mathrm{reg}}$  est existentiel positif dans $A$\fg\hfill{ }
\smallskip

\noindent qui est \'equivalente \`a (C) si $A$ est int\`egre. Nous utiliserons en outre la condition
\smallskip

\noindent(Q)\hfill  \og tout quotient int\`egre de $A$ v\'erifie (C)\fg.\hfill{ }

\begin{prop}\label{EltsReg} Soit $A$ un anneau noeth\'erien. On fait l'une des hypoth\`eses suivantes:
\begin{romlist}
\item\label{EltsReg1} $A$ v\'erifie  \textup{(C)};
\item\label{EltsReg2} pour tout id\'eal premier $\fp$ associ\'e \`a $A$, le quotient $A/\fp$ v\'erifie \textup{(C)}.
\end{romlist}
Alors $A$ v\'erifie $\mathrm{(C^{reg})}$, et v\'erifie  \textup{(C)} s'il est \emph{r\'eduit}.
\end{prop}

\dem soient $\fp_1,\dots,\fp_n$ les id\'eaux premiers associ\'es \`a $A$. Par d\'efinition, chaque $\fp_i$ est l'annulateur d'un \'el\'ement $\alpha_i$ de $A$; on sait en outre (\cite{BourAC3}, chapitre IV, \S1, n\up{\accent'27} 1, cor. 3 de la prop. 2) que $A\setminus A^{\mathrm{reg}}$ est la r\'eunion des $\fp_i$, de sorte que 
$$A^{\mathrm{reg}}=\bigcap_{i=1}^{n}(A\setminus\fp_i).$$
Dans le cas \ref{EltsReg2}, chaque $A\setminus\fp_i$ est existentiel positif par \ref{QuotientDioph}, et il en est de m\^eme dans le cas \ref{EltsReg1} puisque $A\setminus\fp_i=\left\{t\in A\mid t\,\alpha_i\neq0\right\}$. Donc $A$ v\'erifie $\mathrm{(C^{reg})}$. Si $A$ est r\'eduit, il suffit de remarquer que les $\fp_i$ sont les id\'eaux premiers minimaux de $A$ et que $A\setminus\{0\}$ est la r\'eunion des $A\setminus\fp_i$.\qed

\begin{thm}\label{ThNoethQque} Soit $A$ un anneau noeth\'erien v\'erifiant \textup{(Q)}, et soit $S$ une partie multiplicative de $A$. Alors:
\begin{romlist}
\item\label{ThNoethQque1} $A$ v\'erifie \textup{(C)};
\item\label{ThNoethQque2} $S^{-1}A$ v\'erifie \textup{(Q)} \textup{(et donc (C), d'apr\`es \ref{ThNoethQque1})}.
\end{romlist}
\end{thm}

\dem \ref{ThNoethQque1} d'apr\`es \cite{BourAC3} (chap. IV, \S1, n\up{\accent'27} 4, th\'eor\`eme 1) il existe une suite d'id\'eaux
$$A=I_0\supset I_1\supset\cdots\supset I_n=\{0\}$$
et des id\'eaux premiers $\fp_0,\ldots,\fp_{n-1}$ tels que, pour tout $j$, le $A$-module $I_j/I_{j+1}$ soit isomorphe \`a $A/\fp_j$; fixons pour chaque $j$ un \'el\'ement $\alpha_j$ de $I_j$ engendrant $I_j/I_{j+1}$. 
Pour tout $x\in A$, on a alors les \'equivalences
$$
{\renewcommand{\arraystretch}{1.5}
\begin{array}{rcl}
x\neq0 & \Longleftrightarrow & \bigvee_{j=0}^{n-1} \;(x\in I_j\setminus I_{j+1})\\
&  \Longleftrightarrow & \bigvee_{j=0}^{n-1} \;(\exists\,t)\,(\exists\,u) (x=t\alpha_j+u \:\wedge\: t\in A\setminus \fp_j\: \wedge\: u\in I_{j+1})
\end{array}
}$$
d'o\`u la conclusion puisque $A\setminus\fp_j$ et  $I_{j+1}$ sont existentiels positifs, d'apr\`es l'hypoth\`ese sur $A$ et le lemme \ref{QuotientDioph}.\smallskip

\noindent \ref{ThNoethQque2} Les quotients int\`egres de $S^{-1}A$ sont des anneaux de fractions de quotients int\`egres de $A$, donc le cas \ref{CorLocalise2} de \ref{CorLocalise} entra\^{\i}ne que $S^{-1}A$ v\'erifie \textup{(Q)}. \qed\medskip

Rappelons  qu'un anneau $A$ est un \emph{anneau de Jacobson} si tout id\'eal premier de $A$ est intersection d'id\'eaux maximaux. 

\begin{cor}\label{CorJac} Tout anneau de Jacobson noeth\'erien v\'erifie \textup{(Q)}.

Par suite, d'apr\`es \textup{\ref{ThNoethQque}},  tout anneau de fractions d'un anneau de Jacobson noeth\'erien v\'erifie \textup{(C)}. 

En particulier:
\begin{itemize}
\item tout anneau artinien v\'erifie  \textup{(C)};
\item soit $k$  un corps ou un anneau de Dedekind ayant une infinit\'e d'id\'eaux maximaux; alors toute $k$-alg\`ebre essentiellement de type fini v\'erifie \textup{(C)}.
\end{itemize}
\end{cor}

\dem 
Il est clair qu'un anneau de Jacobson int\`egre et local est un corps, de sorte que, d'apr\`es \ref{ThNonHens}, tout anneau de Jacobson \emph{int\`egre} et noeth\'erien v\'erifie (C). Comme il est non moins clair qu'un quotient d'un anneau de Jacobson est de Jacobson, on en d\'eduit la premi\`ere assertion, qui entra\^{\i}ne les autres.\qed

\section{Fonctions analytiques.}

Pour les notions de base sur les espaces analytiques et les espaces de Stein, le lecteur pourra consulter \cite{GrRe79}. Un espace analytique sera toujours suppos\'e  s\'eparable (en particulier l'ensemble de ses composantes irr\'eductibles est d\'enombrable) et de dimension finie.

Parmi les espaces de Stein on compte notamment les espaces $\CC^n$, les polydisques, et les sous-espaces analytiques ferm\'es de ceux-ci.

Si $(X,\cO_X)$ est un espace analytique, nous noterons $\cHol(X)$ l'anneau $H^0(X,\cO_X)$ des fonctions holomorphes globales sur $X$.

\begin{prop}\label{PropExistDiviseur} Soit $(X,\cO_X)$ un espace analytique de Stein, irr\'eductible et r\'eduit, de dimension $>0$, et soit $P$ un point de $X$. Il existe un sous-espace ferm\'e $Y$ de $X$ ayant les propri\'et\'es suivantes:
\begin{romlist}
\item\label{PropExistDiviseur1} $P\not\in Y$, et $Y$ est disjoint du lieu singulier de $X$;
\item\label{PropExistDiviseur2} $Y$ est irr\'eductible et r\'eduit, et son id\'eal $\cI_Y\subset\cO_X$ est localement principal;
\item\label{PropExistDiviseur3} l'id\'eal $\fp=H^0(X,\cI_Y)$ de l'anneau $A:=\cHol(X)$ form\'e des fonctions holomorphes sur $X$ nulles sur $Y$ est de type fini;
\item\label{PropExistDiviseur4} avec les notations de \ref{PropExistDiviseur3}, l'anneau $A_\fp$ est un anneau de valuation discr\`ete.
\end{romlist}
\end{prop}

\dem: soit $Z\subset X$ le sous-espace ferm\'e r\'eunion du lieu singulier de $X$ et du point $P$. Vu les hypoth\`eses sur $X$ il existe un point $Q\in X\setminus Z$. Comme $X$ est de Stein il existe une fonction holomorphe $f$ sur $X$ qui vaut $1$ sur $Z$, $0$ en $Q$ et dont la diff\'erentielle en $Q$ n'est pas nulle. Le sous-espace ferm\'e $Y'$ de $X$ d\'efini par $f=0$ est un diviseur contenu dans $X\setminus Z$, et lisse au voisinage de $Q$. Soit $Y$ l'unique composante irr\'eductible de $Y'$ contenant $Q$. Il est clair que $Y$ v\'erifie les conditions \ref{PropExistDiviseur1}, \ref{PropExistDiviseur2} et \ref{PropExistDiviseur4}, et la propri\'et\'e \ref{PropExistDiviseur3} r\'esulte du lemme \ref{LemSecTF} ci-dessous (appliqu\'e \`a $\cF=\cI_Y$).\qed

\subsubsection{Notations. }\label{NotRangCoh} Soit $(X,\cO_X)$ un espace analytique, et soit $\cF$ un $\cO_X$-module coh\'erent. Pour tout $x\in X$ notons $\cF_x$ le $\cO_{X,x}$-module de type fini fibre de $\cF$ en $x$, et $\cF(x)$ le $\CC$-espace vectoriel de dimension finie $\cF_x/\fm_x\cF_x$ o\`u $\fm_x$ est l'id\'eal maximal de $\cO_{X,x}$. Le lemme de Nakayama implique que $\dim_\CC\cF(x)$ est le plus petit entier $r$ tel qu'il existe un voisinage $U$ de $x$ et une surjection $\cO_U^r\to\cF_{|U}$.

\begin{sublem}\label{LemSecTF} Soit $(X,\cO_X)$ un espace analytique de Stein, et soit $\cF$ un $\cO_X$-module coh\'erent. Pour que le $\cHol(X)$-module $H^0(X,\cF)$ soit de type fini, il faut et il suffit que la fonction $x\mapsto\dim_\CC\cF(x)$ soit born\'ee sur $X$. 
\end{sublem}

\dem la n\'ecessit\'e est \'evidente: si $H^0(X,\cF)$ est engendr\'e par $r$ \'el\'ements, alors la fonction en question est major\'ee par $r$.

Pour montrer la suffisance, on proc\`ede par r\'ecurrence sur la dimension $d$ du support de $\cF$, avec la convention $\dim(\emptyset)=-1$.  Si $d=-1$ alors $\cF$ est nul et il n'y a rien \`a d\'emontrer. 

Soit $r$ un majorant de la fonction $x\mapsto\dim_\CC\cF(x)$. Comme l'ensemble des composantes irr\'eductibles de $\mathrm{Supp}\,(\cF)$ est localement fini, il existe un sous-espace ferm\'e discret  $Z\subset X$ qui rencontre toutes ces composantes. Pour chaque $z\in Z$, soit $(s_{1,z},\ldots,s_{r,z})$ une famille g\'en\'eratrice \`a $r$ \'el\'ements de $\cF(z)$. On obtient donc $r$ sections globales de la restriction de $\cF$ \`a $Z$; comme $X$ est de Stein, elles se rel\`event en $r$ sections $\uple[s]{r}$ de $\cF$ sur $X$, d\'efinissant un morphisme $\varphi:\cO_X^r\to\cF$. Notons $\cF_1$ le conoyau de $\varphi$, qui est un $\cO_X$-module coh\'erent. Par construction, $\varphi$ est surjectif au voisinage de $Z$, de sorte que $\mathrm{Supp}\,(\cF_1)$ ne contient aucune composante de $\mathrm{Supp}\,(\cF)$. Donc $\dim\mathrm{Supp}\,(\cF_1)<d$, et par hypoth\`ese de r\'ecurrence $H^0(X,\cF_1)$ est un $\cHol(X)$-module de type fini (noter que comme $\cF_1$ est un quotient de $\cF$ on a $\dim\cF_1(x)\leq\dim\cF(x)$ pour tout $x\in X$). Comme $X$ est de Stein, la suite exacte $\cO_X^r\to\cF\to\cF_1\to0$ donne une suite exacte $\cHol(X)^r\to H^0(X,\cF)\to H^0(X,\cF_1)\to0$ de $\cHol(X)$-modules, qui montre que $H^0(X,\cF)$ est de type fini. (L'argument montre plus pr\'ecis\'ement qu'il est engendr\'e par $r(d+1)$ \'el\'ements).\qed

\begin{thm}\label{ThStein} Soit $(X,\cO_X)$ un espace analytique de Stein, irr\'eductible et r\'eduit. Alors l'anneau $\cHol(X)$ v\'erifie \textup{(C)}.
\end{thm}

\dem posons $A=\cHol(X)$, et proc\'edons par  r\'ecurrence sur $d:=\dim X$. Si $d=0$, alors $A=\CC$. Si $d>0$, soient $P\in X$ et $Y\subset X$ comme dans la proposition \ref{PropExistDiviseur}. On note $\fp$ l'id\'eal de $Y$ dans $A$, et $\fm$ l'id\'eal maximal du point $P$. Appliquons le lemme \ref{LemDeuxId} avec $\fp_1=\fp$ et $\fp_2=\fm$. La proposition  \ref{PropExistDiviseur} assure que les conditions \ref{LemDeuxId1} et \ref{LemDeuxId2} de \ref{LemDeuxId} sont satisfaites (cf. la remarque \ref{RemLemDeuxId}\,(2)). 
De plus $\fp_1$ est de type fini d'apr\`es \ref{PropExistDiviseur}\,\ref{PropExistDiviseur3}, et il en est de m\^eme de $\fp_2$ d'apr\`es le lemme \ref{LemSecTF}. L'anneau $A/\fp_1$ s'identifie \`a $\cHol(Y)$ donc v\'erifie (C) par hypoth\`ese de r\'ecurrence, et $A/\fp_2\cong\CC$ v\'erifie trivialement (C), donc  le lemme \ref{LemDeuxId} donne bien le r\'esultat voulu.\qed
\bibliography{biblio}
\end{document}